\documentclass[twoside,11pt]{article}

\usepackage{jmlr2earxiv}

\usepackage{amsmath}
\usepackage{pgfplots}
\usepackage{subcaption}
\usepackage{graphicx}

\usepackage[noend]{algpseudocode}
\usepackage{algorithm}
\usepackage{algorithmicx}
\newcommand{\lt}{\left}
\newcommand{\rt}{\right}
\newcommand*\Let[2]{\State #1 $\gets$ #2}
\algrenewcommand\algorithmicrequire{\textbf{Input:}}
\algrenewcommand\algorithmicensure{\textbf{Output:}}

\def\balpha{{\boldsymbol \alpha}}
\def\bPhi{{\boldsymbol \Phi}}
\def\brho{{\boldsymbol \rho}}

\def\bEps{{\boldsymbol \epsilon}}
\def\bB{{\mathbf B}}
\def\bb{{\mathbf b}}
\def\bX{{\mathbf X}}
\def\bx{{\mathbf x}}
\def\bone{{\mathbf 1}}
\def\bt{{\mathbf t}}
\def\bw{{\mathbf w}}
\def\bi{{\mathbf i}}

\def\C{{\mathbb C}}
\DeclareMathOperator*{\argmin}{argmin}
\DeclareMathOperator{\diag}{diag}
\DeclareMathOperator{\Diag}{Diag}
\def\opt{{\mbox{opt}}}

\def\R{{\mathbb R}}

\jmlrheading{1}{2017}{1--1}{11/17}{??/18}{askham17a}{Travis Askham, Peng Zheng, Aleksandr Aravkin, and J. Nathan Kutz}

\ShortHeadings{Robust and scalable DMD}{Askham et al.}
\firstpageno{1}

\begin{document}

\title{Robust and scalable methods for the dynamic mode
decomposition}

\author{\name Travis Askham \email askham@uw.edu \\
       \addr Department of Applied Mathematics\\
       University of Washington\\
       Seattle, WA 98195-3925, USA
       \AND
       \name Peng Zheng \email zhengp@uw.edu \\
       \addr Department of Applied Mathematics\\
       University of Washington\\
       Seattle, WA 98195-3925, USA
       \AND
       \name Aleksandr Aravkin \email saravkin@uw.edu \\
       \addr Department of Applied Mathematics\\
       University of Washington\\
       Seattle, WA 98195-3925, USA
       \AND
       \name J. Nathan Kutz \email kutz@uw.edu \\
       \addr Department of Applied Mathematics\\
       University of Washington\\
       Seattle, WA 98195-3925, USA}

\editor{ }

\maketitle

\begin{abstract}
  The dynamic mode decomposition (DMD) is a broadly 
  applicable dimensionality reduction algorithm that 
  approximates a matrix containing time-series data by the outer
  product of a matrix of exponentials, representing
  Fourier-like time dynamics, and a matrix of coefficients,
  representing spatial structures. This interpretable 
  spatio-temporal decomposition is 
  commonly computed using linear algebraic techniques in its
  simplest formulation or a nonlinear optimization procedure
  within the variable projection framework. For data with sparse
  outliers or data which are not well-represented by exponentials
  in time, the standard Frobenius norm fit of the data creates
  significant biases in the recovered time dynamics.
  As a result, practitioners are left to clean such defects
  from the data manually or to use a black-box cleaning approach
  like robust PCA. As an alternative, we propose a framework 
  and a set of algorithms for 
  incorporating robust features into the nonlinear 
  optimization used to compute the DMD itself. The algorithms 
  presented are flexible, allowing for regularizers and 
  constraints on the optimization, and scalable, using a 
  stochastic approach to decrease the computational cost 
  for data in high dimensional space. Both synthetic and real
  data examples are provided.
\end{abstract}

\begin{keywords}
  robust statistics, dynamic mode decomposition,
  scalable algorithms
\end{keywords}

\section{Introduction}

Dimensionality reduction is a critically enabling tool
in machine learning applications.
Specifically, extracting the dominant
low-rank features from a high-dimensional data matrix
${\mathbf X} \in \R^{m\times n}$ allows one to efficiently
perform tasks associated
with clustering, classification and prediction. As defined by
\cite{cunningham2015linear}, {\em linear} dimensionality
reduction methods solve an optimization problem with objective
$f_{{\mathbf X}}(\cdot)$ over a manifold ${\mathcal M}$ to produce a
linear transformation $P$ which maps the columns of
${\mathbf X}$ to a lower dimensional space. 
Many popular methods can be written in this framework
by an appropriate definition of $f_{{\mathbf X}}(\cdot)$ and
specification of the manifold $\mathcal{M}$. For instance,
the principal component analysis (PCA) may be written
as

\begin{equation} \label{eq:lindr}
  \hat{M} = \argmin_{M \in \mathcal{M}} \| {\mathbf X}
  - MM^\intercal {\mathbf X} \|_F \, , \qquad \mathcal{M} =
  \mathcal{O}^{m\times k} \; ,
\end{equation}
where $\mathcal{O}^{m\times k}$ is the manifold of
$m\times k$ matrices with orthonormal columns,
i.e. $\mathcal{M}$ is a
Stiefel manifold. The map $P$ is then given by $\hat{M}^\intercal$. 
One of the primary
conclusions of the survey \citep{cunningham2015linear},
is that --- aside from the PCA itself --- many of the
common methods for linear dimensionality
reduction based on eigenvalue solvers are actually
sub-optimal heuristics and the direct solution of the
optimization problem \eqref{eq:lindr} should be
preferred. 

In this manuscript, we consider a particular linear
dimensionality reduction technique: the dynamic
mode decomposition (DMD). In the past
decade, the DMD has been applied to the analysis of
fluid flow experiments and simulations, machine learning
enabled control systems, and Koopman spectral analysis,
among other data-intensive problems described by dynamical
systems. The success of the algorithm is largely due 
to the interpretability of the low-rank spatio-temporal 
modes it generates in approximating the dominant features 
of the data matrix ${\bf X}$.
The DMD was originally defined to be the
output of an algorithm for characterizing time-series
measurements of fluid flow
data~\citep{Schmid2008aps,schmid2010dynamic}. It was
later reformulated by \cite{tu2014dynamic} as a
least-squares regression problem whereby the DMD could
be stably computed using a Moore-Penrose pseudo-inverse
and an eigenvalue decomposition. An earlier though less
commonly used formulation, the {\em optimized}
DMD \citep{chen2012}, can be phrased as the optimization
problem
\begin{equation} \label{eq:optdmdintro}
\hat{M} = \argmin_{M \in \mathcal{M}} \| {\mathbf X}
  - MM^\dagger {\mathbf X} \|_F \; , \quad
  \mathcal{M} = {\boldsymbol \Phi}(\C^k) \; ,
\end{equation}
where the map $\balpha \mapsto {\bf \Phi}(\balpha)$
defines a matrix with columns corresponding
to exponential time dynamics (see Section~\ref{sec:dmd})
and $M^\dagger$ denotes the Moore-Penrose pseudo-inverse
of $M$.
This can be thought of as a best-fit linear dynamical
system approximation of the data. In agreement with the
conclusions of \cite{cunningham2015linear},
the optimized DMD, while more costly to compute,
is more robust to additive noise
than the exact DMD and its noise-corrected alternatives
\citep{dawson2016,askham2017variable}. It is also more flexible than
the exact DMD, allowing for non-equispaced snapshots.
While the
optimized DMD does not fit directly into the
optimization framework of \cite{cunningham2015linear},
which is defined for $\mathcal{M}$ either a Stiefel
manifold or a Grassmannian manifold, it can be computed
efficiently using classical variable projection methods
\citep{golub1979,askham2017variable}.

The DMD has been used in a variety of fields
where the nature of the data can lead to corrupt
and noisy measurements.  This includes applications ranging
from neuroscience~\citep{brunton2016extracting} to video
processing~\citep{jake,erichson2015compressed} to fluid
dynamics~\citep{Schmid2008aps,schmid2010dynamic,
  gueniat2015a,dawson2016}.
Although the Frobenius norm used in the definition of
the optimized DMD \eqref{eq:optdmdintro} is appealing
due to its physical interpretability in many applications
(energy, mass, etc.), it has significant flaws that
can severely limit its applicability.  Specifically, corrupt
data or large noise fluctuations can lead to significant
deformation of the DMD approximation of the data
because the Frobenius norm implicitly assigns a very low
probability to such outliers
(see Section~\ref{sec:robuststats}). In practice, these
outliers are often removed from the data
manually or using a black-box filtering approach like
robust PCA
\citep{locantore1999robust,wright2009robust,candes2011robust}.
However, such approaches ignore the structure
of the DMD approximation and may introduce biases
of their own. Further, it is desirable that DMD methods
not only be robust to ``noisy'' outliers but also to
non-exponential structure in the data. We therefore
develop an alternative approach to increase the
robustness of the DMD. In particular, 
we modify the optimized DMD definition \eqref{eq:optdmdintro}
to incorporate ideas from the field of robust
statistics~\citep{maronna2006robust,huber2011robust}
in order to produce a decomposition that is significantly
less sensitive to outliers in the data.

Because the new problem formulation incorporates
robust norms, many of the efficient strategies
used in variable projection algorithms for problems
defined in the Frobenius norm are no longer available.
To remedy this, we develop a number of algorithms
based on modern variable projection
methods~\citep{aravkin2012,aravkin2017efficient}
which exploit the structure of the DMD for increased
performance. In particular, we can incorporate nonsmooth features,
such as regularizers and constraints, 
and scale to large problems using stochastic variance reduction techniques. 

This flexible architecture
allows us to impose physically relevant constraints
on the optimization that are critical for tasks such
as future-state prediction. For instance, we can
impose the constraint that the real parts of the DMD
eigenvalues are non-positive, thus ensuring that
solutions do not grow to infinity when forecasting.

The effect of noise on the DMD is a well-studied
area. Controlling for the bias of the exact DMD
in the presence of additive noise was treated by
\cite{hemati2017} and \cite{dawson2016}. A Bayesian
formulation of the DMD was presented by
\cite{takeishi2017bayesian}. This formulation is
flexible enough to incorporate robust statistics
but this was not a focus of that work.
\cite{dicle2016robust} presented a robust formulation
of exact DMD type, which complements the current work.

The rest of this manuscript is organized as follows.
In Section~\ref{sec:prelim}, we provide some necessary
preliminaries from the DMD, robust statistics, and variable
projection literature and we present our problem formulation.
A detailed description of the algorithms we use to solve the robust DMD formulation 
follows in Section~\ref{sec:methods}. We apply these methods to
synthetic data in Section~\ref{sec:synthetic} and to real
data in Section~\ref{sec:real}. Finally, we provide some
concluding remarks and describe possible future directions
in Section~\ref{sec:conclusion}.

\section{Preliminaries}
\label{sec:prelim}
In this section, we outline some of the precursors of
this work and present our problem formulation.

\subsection{Dynamic mode decomposition}
\label{sec:dmd}
As mentioned above, the dynamic mode decomposition (DMD)
corresponds to a best-fit linear dynamical model of the
data. Because linear dynamics produce exponential functions
in time, the DMD may be written as an exponential fitting problem. 
Let $\bX \in \C^{m\times n}$ be a snapshot matrix whose rows 
$\bx_{j}$ are samples of an $n$ dimensional dynamical 
system at a set of $m$ sample times $t_j$. For a given rank 
$k$, let $\balpha \in \C^k$ be a vector of complex numbers
specifying time dynamics. We then define the matrix 
$\bPhi(\balpha ; \bt)$ by 

\begin{equation} \label{eq:phidef}
  \Phi_{ij} (\balpha) = e^{\alpha_j t_i} \; .
\end{equation}
When it is clear in context, we often drop the dependence
of $\bPhi$ on $\balpha$ and $\bt$. 

Let $\bB \in \C^{k \times n}$ be a matrix of coefficients for 
the exponential fit. The so-called \emph{optimized 
DMD} \citep[see][]{chen2012} is defined to be the solution 
of the following optimization problem:
\begin{equation} \label{eq:optdmd}
\displaystyle \min_{\balpha, \bB} \dfrac12 \| \bX - \bPhi(\balpha) \bB \|_F^2 \; .
\end{equation}
The problem \eqref{eq:optdmd} is a large, nonlinear least squares problem; 
in particular it is highly non-convex. 
The classical variable projection framework provides an efficient method
for computing a (local) solution. 

Let 
\begin{equation}
f_\opt (\balpha,\bB) = 
\dfrac12 \| \bX - \bPhi(\balpha) \bB \|_F^2 \; . \nonumber
\end{equation}
The classical variable projection framework is based on 
the observation that for a fixed $\balpha$, it is easy 
to optimize $f_\opt$ in $\bB$. In fact, for the least squares case, 
we have a closed form expression 
 \begin{equation}
 \label{eq:ls-vp}
 \bB(\balpha) := \arg\min_{\bB} f_\opt(\balpha, \bB) =  \bPhi(\balpha)^\dagger \bX,
 \end{equation}
 where  $\bPhi(\balpha)^\dagger$ denotes the Moore-Penrose 
pseudo-inverse of $\bPhi(\balpha)$. Let 
\begin{equation}
  \tilde{f}_\opt (\balpha) = \min_{\bB} f_\opt(\balpha, \bB) :=
  \dfrac12 \| \bX - \bPhi(\balpha)  \bB(\balpha)  \|_F^2.
  \nonumber
\end{equation}
The variable projection (VP) technique finds the
minimizer of $\tilde{f}_\opt(\balpha)$ using an iterative
method. First and second derivatives of $\tilde{f}$ with respect to $\balpha$ are easily computed
\citep[see Theorem 2 of][]{bell2008algorithmic}: 
\begin{equation}
\label{eq:deriv_form}
\begin{aligned}
  \nabla_{\balpha} \tilde{f}_\opt (\balpha) & =   \partial_{\balpha} f_\opt \vert_{\balpha,\bB(\balpha)} \\
  \nabla_{\balpha}^2 \tilde{f}_\opt(\balpha) & =  \lt.\left[\partial^2_{\balpha}f_\opt - \partial_{\balpha, \bB}f_\opt(\partial^2_{\bB}f_\opt)^{-1}\partial_{\bB,\balpha}f_\opt\right]\rt\vert_{\balpha,\bB(\balpha)}.
\end{aligned}
\end{equation}
These formulas allow first- and second-order methods to be directly applied to $\tilde{f}_\opt$, 
including steepest descent, BFGS, and Newton's method. 
The matrix $\bB(\balpha)$ is updated every time $\balpha$ changes.  
Gauss-Newton and Levenberg-Marquardt (LM) have been classically used for exponential fitting; these methods do not use the Hessian in~\eqref{eq:deriv_form}, opting for simpler approximations. The method was used for exponential fitting by~\cite{golub1979}.
 While VP originally referred to 
least-squares projection (using the closed-form solution $\bB(\balpha)$ in~\eqref{eq:ls-vp}), follow-up work 
considered more general loss functions, using the term {\it projection} to refer to partial minimization~\citep{aravkin2012,aravkin2017efficient}.

For practitioners, the optimized DMD may be less familiar
than exact DMD~\citep{tu2014dynamic}. 
We favor the optimized DMD for its
performance on data with additive noise 
\citep[see][]{askham2017variable} and its flexibility.  In particular, the optimized formulation enables the 
contributions of the current work.  
For a review of the DMD and its applications, see \citet{tu2014dynamic} and \cite{kutz2016}.

\subsection{Robust Formulations}
\label{sec:robuststats}

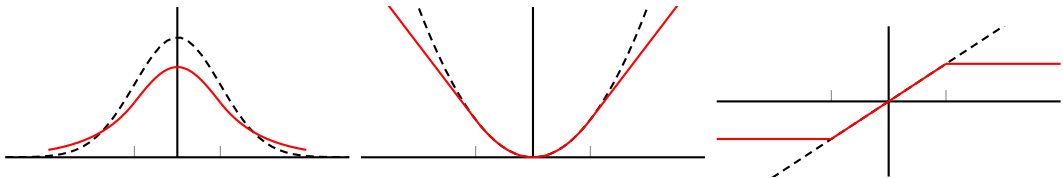
\begin{figure} \label{GLT-KF}
\centering
\begin{tikzpicture}
  \begin{axis}[
    thick,
    width=.3\textwidth, height=2cm,
    xmin=-4,xmax=4,ymin=0,ymax=0.5,
    no markers,
    samples=50,
    axis lines*=left, 
    axis lines*=middle, 
    scale only axis,
    xtick={-1,1},
    xticklabels={},
    ytick={0},
    ] 
\addplot[domain=-4:+4,densely dashed]{exp(-.5*x^2)/sqrt(2*pi)};
\addplot[red, domain=-3:-1]{0.3*exp(-abs(x)+.5)};
  \addplot[red, domain=-1:+1]{0.3*exp(-.5*x^2)};
 \addplot[red, domain=1:+3]{0.3*exp(-abs(x)+.5)};

  \end{axis}
\end{tikzpicture}
\begin{tikzpicture}
  \begin{axis}[
    thick,
    width=.3\textwidth, height=2cm,
    xmin=-3,xmax=3,ymin=0,ymax=2,
    no markers,
    samples=50,
    axis lines*=left, 
    axis lines*=middle, 
    scale only axis,
    xtick={-1,1},
    xticklabels={},
    ytick={0},
    ] 
\addplot[domain=-3:+3,densely dashed]{.5*x^2};
 \addplot[red, domain=-3:-1]{abs(x)-.5};
  \addplot[red, domain=-1:+1]{.5*x^2};
 \addplot[red, domain=1:+3]{abs(x)-.5};

  \end{axis}
\end{tikzpicture}
\begin{tikzpicture}
  \begin{axis}[
    thick,
    width=.3\textwidth, height=2cm,
    xmin=-3,xmax=3,ymin=-2,ymax=2,
    no markers,
    samples=50,
    axis lines*=left, 
    axis lines*=middle, 
    scale only axis,
    xtick={-1,1},
    xticklabels={},
    ytick={0},
    ] 
\addplot[domain=-3:3,densely dashed]{x};
 \addplot[red, domain=-3:-1]{-1};
  \addplot[red, domain=-1:+1]{x};
 \addplot[red, domain=1:+3]{1};
  \end{axis}
\end{tikzpicture}
    \caption{\label{fig:robust}Gaussian (black dash) and Huber (red solid) Densities, Negative Log Likelihoods, and Influence Functions.}
\end{figure}

The optimized DMD problem~\eqref{eq:optdmd} is formulated using the least-squares error norm, 
which is equivalent to assuming a Gaussian model on the errors between predicted and observed data: 
\[
\bX = \bPhi(\balpha)\bB + \bEps, \quad \bEps \sim N(0, \sigma^2 I). 
\]
This error model, and the corresponding formulation, are vulnerable to outliers in the data. 
Both DMD and optimized DMD are known to be sensitive to outliers, so in practice data are `pre-cleaned' 
before applying these approaches. 

In many domains, formulations based on robust statics have become the method of choice for dealing with contaminated 
data. Two common approaches are 
\begin{itemize}
\item to replace the LS penalty with one that penalizes deviations 
  less harshly and
\item to solve an extended problem that explicitly identifies 
outliers while fitting the model.
\end{itemize} 
The first approach, often called M-estimation~\citep{huber2011robust,maronna2006robust}, is illustrated in Figure~\ref{fig:robust}. Replacing the least squares penalty by the Huber penalty
\[
\rho(z) = \begin{cases} 
\frac{1}{2} |z|^2 & \mbox{if} \; |z| \leq \kappa \\
\kappa |z| - \frac{1}{2}\kappa^2 & \mbox{if} \; |r| > \kappa
\end{cases}
\]
corresponds to choosing the solid red penalty rather than the dotted black least squares penalty in the center panel of Figure~\ref{fig:robust}.
This corresponds to modeling errors $\bEps$ using the density $\exp(-\rho)$, which has heavier tails than the Gaussian (see left panel of Figure~\ref{fig:robust}). Heavier tails means deviations (i.e. larger residuals) are more likely than under the Gaussian model, and so observations
that deviate from the norm have less {\it influence}, i.e. effect on the fitted parameters $(\balpha, \bB)$ than under the Gaussian model (see right panel of Figure~\ref{fig:robust}). The M-estimator-DMD problem can be written as 
\[
\min_{\balpha, \bB} \sum_{j=1}^n \rho(X_{\cdot j} - \bPhi(\balpha)\bB_j) := \sum_{j=1}^n\rho_j(\balpha, \bB),
\]
where the sum is run across columns. 

Another approach, called trimmed estimation, builds on M-estimation by coupling explicit outlier identification/removal with model fitting. 
The trimmed DMD formulation for any penalty $\rho$ is given by 
\begin{align}\label{eq:trimmed}
\min_{\balpha, \bB} \sum_{l=1}^h \rho_{j_l}(\balpha, \bB),
\end{align}
where $\rho_{j_1}(\balpha, \bB) \leq \cdots \leq \rho_{j_h}(\balpha, \bB)$ are the first $h$ order statistics of the objective values and $\{j_1, \ldots, j_h\}\subseteq\{1,\ldots,n\}$. Interpreting the loss $\rho_j$ as the negative log likelihood of the $j$th observed column, it is clear that trimming  jointly fits a likelihood model while simultaneously eliminating the influence of all low-likelihood observations. An equivalent formulation to~\eqref{eq:trimmed} replaces the order statistics with explicit weights 
\begin{align}\label{eq:trimmedw}
\min_{\balpha, \bB, \bw} \sum_{j=1}^n w_j \rho_j(\balpha, \bB), \quad 0 \leq w_j \leq 1, \quad \bone^\intercal\bw = h.
\end{align}
The reader should verify that~\eqref{eq:trimmedw} and~\eqref{eq:trimmed} are equivalent.

Trimmed M-estimators were initially introduced by~\cite{rousseeuw1985multivariate} in the context of least-squares regression. The author's original motivation was to develop linear regression estimators that have a high breakdown point (in this case 50\%) and good statistical efficiency (in this case $n^{-1/2}$)\footnote{Breakdown refers to the percentage of outlying points which can be added to a dataset before the resulting M-estimator can change in an unbounded way.}. For a number of years, the difficulty of efficiently optimizing LTS problems limited their application. 
However, recent work has made it possible to efficiently apply trimming to general models~\citep{yang2016high,aravkin2016smart}. We 
show how to incorporate trimming into the robust DMD framework.

\subsection{Regularization}
\label{sec:reg}
 Optimized DMD allows prior knowledge to be incorporated into the optimization formulation, 
either through constraints on variables, or regularization terms. In all exponential fitting problems, the real parts of $\balpha$ coefficients
play a major role in explaining the data because of the exponential growth of $\bPhi(\balpha)$. 
A natural regularization is to restrict the magnitudes of the real parts of $\balpha$, imposing the constraint $\mbox{real}(\balpha) \leq \gamma$
with $\gamma$ chosen by the user.  We write the constraint as follows: 
\[
r(\balpha) = 
  \begin{cases}
   0 & \mbox{if} \; \mbox{real}(\balpha) \leq \gamma\\
   \infty & \mbox{if} \; \mbox{real}(\balpha) > \gamma.
   \end{cases}
\] 
This is a simple convex function that admits a trivial proximal operator (see~\cite{combettes2011proximal}):
the projection onto the shifted left half-plane in $\mathbb{C}^k$. The VP technique can be easily adapted to incorporate such functions 
on $\balpha$. 

Constraints and penalties can also be imposed on the matrix $\bB$. 
We assume that only smooth separable regularization penalties can be used; 
and in this case, the regularization is added to the $g$ function. 

\subsection{Problem formulation}
\label{sec:problemform}
Let $h(\bB)$ and $r(\balpha)$ be
convex regularization terms. We formulate the general robust DMD problem as follows:

\begin{equation} \label{eq:problemform}
\min_{\balpha, \bB, \bw}  f(\balpha,\bB, \bw) :=  g(\balpha,\bB, \bw) + r(\balpha) + s(\bw)\; ,
\end{equation}
where $r(\balpha)$ encodes optional regularization functions (or constraints) 
 for $\balpha$ (see Section~\ref{sec:reg}) and 
\begin{equation} \label{eq:formg}
  g(\balpha,\bB, \bw) = \sum_{j=1}^n w_j\rho\left ( X_{\cdot j} 
  - \bPhi(\balpha) \bB_{\cdot j} \right) + q(\bB_{\cdot j})
\end{equation}
with  $\rho$ any differentiable penalty, $q(\bB_{\cdot j})$ representing potential regularizer for columns of $\bB$, and $s(\bw)$ encoding the capped simplex constraints: 
   \begin{equation}
   \label{eq:simplex}
   s(\bw) = \begin{cases}
   0 & \mbox{if} \; 0 \leq w_j \leq 1, \; \bone^\intercal\bw = h \\
   \infty & \mbox{else}.
   \end{cases}
   \end{equation}
These constraints are explained in Section~\ref{sec:robuststats}. The $\bw$ variables select the best-fit 
$h$ columns of the data, and only use those values to update $\balpha$. Since each $w_j \in [0,1]$ rather than 
\{0,1\}, the solutions do not have to be integral. However, for any fixed $(\bB, \balpha)$ there exists a vertex 
solution, since the subproblem in $\bw$ with the other variables fixed is a linear program. 
The function $s(\bw)$ admits a simple proximal operator, which is the projection onto the intersection 
of the $h$-simplex with the unit cube\footnote{This set is called the {\it capped simplex}, and admits fast projections~\citep{aravkin2016smart}.}. 

Setting $h=n$ forces $w_j = 1$ for each column, eliminating 
trimming completely, and reducing~\eqref{eq:problemform}  to a simpler regularized M-estimation form of DMD. 
  
    

For notational convenience, we define a matrix-valued penalty function 
\[
\brho(A) := \begin{bmatrix} 
\rho(A_{1,1}) & \cdots & \rho(A_{1,n}) \\
\vdots & \ddots & \vdots \\
\rho(A_{m,1}) & \cdots & \rho(A_{m,n})
\end{bmatrix}.
\]
In this notation, we can write 
\[
g(\balpha, \bB, \bw) = \bone^\intercal \brho(\bX - \bPhi(\balpha)\bB)\bw + q(\bB),
\]
which makes derivative computations straightforward. 


Our numerical examples use constraints for $\balpha$, but do not regularize $\bB$, that is, $q(\bB) \equiv 0$. 
 However, we consider separable penalties $q$ in the algorithmic description to preserve the generality of~\eqref{eq:problemform}.
 
 \subsection{Gradient computations}
\label{sec:gradient}
We need to compute the gradient of the penalty function \eqref{eq:formg}
with respect to the entries of $\balpha$ and 
$\bB$. In all methods, we treat treat the real
and imaginary components of $\alpha_j$ and $B_{ji}$ as independent, real-valued parameters.

Consider a complex number $z = x + \bi y$.
We write derivative formulas in the Wirtinger sense, 
computing partial derivatives with respect to 
the complex variables. The derivatives for the real 
components can then be recovered from the formulas
\begin{equation}
\dfrac{\partial}{\partial z} = 
\dfrac12 \left ( \dfrac{\partial}{\partial x}
- \bi \dfrac{\partial}{\partial y} \right ).
\end{equation}
%
Let $g(z)$ be a 
function of $z$ which can be written as $g(z) = G(z,\bar{z})$
where $G$ is differentiable with respect to both $z$ 
and $\bar{z}$. The Wirtinger derivative of $g$ is then
the partial derivative of $G$ with respect to $z$, treating
$\bar{z}$ as a constant. For example, the Huber penalty
may be written as 
\[
  \rho(z) = H(z, \bar z; \kappa) = \begin{cases}
   \kappa \sqrt{z\bar{z}} - \frac{1}{2} \kappa^2 
      ,&  |z| \geq \kappa \\
      \frac{1}{2} z\bar{z} ,& |z| < \kappa
    \end{cases} \; .
\]
The Wirtinger derivative of the Huber penalty is then

\[
  \rho'(z) = \dfrac{\partial}{\partial z} H(z, \bar z; \kappa) = 
    \begin{cases} \dfrac{\kappa\bar{z}}{2\sqrt{z\bar{z}}}
      ,& |z| < \kappa \\
      \frac{1}{2} \bar{z}, & |z| \geq \kappa
    \end{cases} \; .
\]
The gradients of $f$ with respect to $\balpha$ and $\bB$ can then be computed using the chain rule: 

%
\begin{equation}
\label{eq:rho_grad}
\begin{aligned}
\nabla_\balpha g(\balpha, \bB, \bw) &= -\diag\left[\bB\Diag(\bw)\brho'(\bX - \bPhi\bB)^\intercal(\Diag(\bt)\bPhi)\right]\\
\nabla_\bB g(\balpha, \bB, \bw) &= -\bPhi^\intercal\brho'(\bX - \bPhi\bB)\Diag(\bw) + \bB^\intercal\nabla q(\bB)\\
\nabla_\bw g(\balpha, \bB, \bw) & =\brho(\bX - \bPhi\bB)^\intercal\bone,
\end{aligned}
\end{equation}
where we define 
\[
\begin{aligned}
\text{diag}(A) & := \begin{bmatrix} a_{11} \\ \vdots \\ \vdots \\ a_{nn}\end{bmatrix},  \qquad
 \text{Diag}(a)  := \begin{bmatrix} a_1 & 0 & \dots & 0 \\
 0& a_2 & \ddots & \vdots \\
 \vdots & \ddots & \ddots& 0 \\
 0 & \cdots & 0 & a_n\end{bmatrix}.
 \end{aligned}
\]
 
\section{Methods}

\label{sec:methods}
In this section, we develop numerical approaches for \eqref{eq:problemform}. 

\subsection{Variable projection framework}

To compute the robust optimized DMD, we apply the variable
projection (VP) technique to the optimization
problem \eqref{eq:problemform}.  Define the reduced function $\tilde f$ and implicit solution $\bB(\balpha)$ by
\begin{equation} \label{eq:reduced}
\begin{aligned}
  \tilde{f}(\balpha, \bw) &= \min_\bB f(\balpha,\bB, \bw) \; ,\\
  \bB(\balpha, \bw) &= \argmin_\bB f(\balpha, \bB, \bw) \; ,
\end{aligned}
\end{equation}
where $f$ is as defined in \eqref{eq:problemform}. The gradient formula~\eqref{eq:deriv_form} holds 
for a very broad problem class. In particular, it holds as long as the following conditions are satisfied~\cite[Theorem 10.58]{rockafellar2009variational}: 
\begin{enumerate}
\item $g(\balpha, \bB, \bw)$ is level-bounded in $\bB$ locally uniformly in $\balpha$; in particular for any 
compact subset of $\balpha$, the union of sublevel sets $\{\bB: g(\balpha, \bB, \bw) \leq \gamma\}$ is bounded. 
\item The gradient of $g(\balpha, \bB, \bw) $ exists and is continuous for all $(\balpha, \bB, \bw)$.  
\item $\bB(\balpha, \bw)$ is unique. 
\end{enumerate}
Several assumptions on $g$, $\bPhi$, and $q$ can be made to ensure these conditions hold. For example, if 
$g$ is differentiable, convex, coercive\footnote{A function $g$ is coercive if it grows in every direction, i.e. 
 $\displaystyle\lim_{\alpha \uparrow \infty} g(\alpha x) = \infty$ for any $x \neq 0$.} in $\bB$, and $\bPhi(\balpha)$ has full rank, then the result holds. If the same conditions hold for $g$, $\bPhi(\balpha)$ does not have full rank, but $q$ is strictly convex, the result holds as well. 
The derivative formulas are valid for all of the examples in the paper, and we have 
\begin{equation}
  \nabla \tilde{f}(\balpha, \bw) =
  \partial_{\balpha,\bw} f(\balpha,\bB, \bw) \vert_{\balpha,\bB(\balpha), \bw} \; .
\end{equation}

We refer to partially minimizing over $f$ over $\bB$ as the {\it inner problem} and minimizing $\tilde f$ as the {\it outer problem}.
When $f$ is convex and smooth with respect to $\bB$, a lot of fast optimization algorithms can be applied to the inner problem.
The inner problem is also embarrassingly parallelizable, and we make use of the problem structure in algorithm design. 
The general VP strategy is to use an iterative method to compute
a (local) minimizer of the reduced function~\eqref{eq:reduced}. 

Solving \eqref{eq:problemform} requires optimization procedures for both the inner and outer problems. 
We discuss these algorithms in the next two subsections.
\subsection{Batch methods}
Consider the problem class where $g$ in \eqref{eq:problemform} is convex and continuously differentiable with respect to $\bB$.
In this case, the inner problem decouples into $n$ independent subproblems:
\begin{equation}
\label{eq:sub}
\bb_j(\balpha, \bw) = \argmin_\bb \quad w_j\rho\lt(X_{\cdot j} - \bPhi(\balpha) \bb\rt) + q(\bb), \quad j = 1,\ldots,n.
\end{equation}
We use BFGS to solve each of these subproblems, since the dimension of each problem is relatively small, 
and BFGS gives a superlinear convergence rate while using only gradient information. 
When $r$ in \eqref{eq:problemform} is continuously differentiable, we 
can also use BFGS as our solver, see Algorithm~\ref{alg:BFGS-BFGS}.
When $r$ is non-smooth but admits an efficient prox operator, 
a first order method such as the proximal gradient method or its accelerations, such as FISTA~\citep{beck2009fast}, 
can be used instead, see Algorithm~\ref{alg:FISTA-BFGS}. We let $\nu$ denote the iteration counter. 

\begin{algorithm}[H]
  \caption[Caption]{\label{alg:BFGS-BFGS}VP using BFGS for outer problem (smooth $r$).}
  \begin{algorithmic}[1]
   \Require{$\balpha^0$, $\bB^0$, $\bw^0$, $H_\balpha^0 = I$, $\nu=0$.}
   \While{not converged}
   	\For{$j = 1,\ldots,n$}
   		\Let{$\bb^{\nu+1}_j$}{$\displaystyle\arg\min_\bb w_j^\nu \rho\lt(X_{\cdot j} - \bPhi(\balpha^\nu) \bb\rt) + q(\bb)$ }
	\EndFor
	\Let{$\bw^{\nu+1}$}{weights update}
	\Let{$f_\balpha^\nu$}{$f(\balpha^\nu, \bB^{\nu+1}, \bw^{\nu+1})$}
	\Let{$g_\balpha^{\nu}$}{$\nabla_\balpha f(\balpha^\nu, \bB^{\nu+1}, \bw^{\nu+1})$}
	\If{$\nu \ge 1$}
	\Let{$s^\nu$}{$f_\balpha^\nu - f_\balpha^{\nu-1}$}
	\Let{$y^\nu$}{$g_\balpha^\nu - g_\balpha^{\nu-1}$}
	\Let{$\beta^\nu$}{$(\lt\langle s^\nu, y^\nu\rt\rangle)^{-1}$}
	\Let{$H_\balpha^\nu$}{$\lt[I - \beta^\nu (s^\nu)(y^\nu)^\intercal\rt]H^{\nu-1}\lt[I - \beta^\nu (y^\nu)(s^\nu)^\intercal\rt]+ \beta (s^\nu)(s^\nu)^\intercal$}
	\EndIf
	\Let{$\balpha^{\nu+1}$}{LineSearch$\lt(\balpha^\nu - \eta_\balpha H_\balpha^\nu g_\balpha^{\nu}\rt)$}
	\Let{$\nu$}{$\nu+1$} 
   \EndWhile
   \Ensure{$\balpha^\nu$, $\bB^\nu$.}
  \end{algorithmic}
\end{algorithm}

\begin{algorithm}[H]
  \caption[Caption]{\label{alg:FISTA-BFGS} VP using prox-gradient for outer problem (prox-friendly $r$).}
  \begin{algorithmic}[1]
   \Require{$\balpha^0$, $\bB^0$, $\bw^0$, $\nu=0$.}
   \While{not converged}
   	\For{$j = 1,\ldots,n$}
   		\Let{$\bb^{\nu+1}_j$}{$\displaystyle\arg\min_\bb w_j^\nu \rho\lt(X_{\cdot j} - \bPhi(\balpha^\nu) \bb\rt) + q(\bb)$ }
	\EndFor
	\Let{$\bw^{\nu+1}$}{weights update}
	\Let{$\balpha^{\nu+1}$}{$\text{prox}_{\eta_\balpha r}\lt(\balpha^\nu - \eta_\balpha \nabla_\balpha f(\balpha^\nu, \bB^{\nu+1}, \bw^{\nu+1})\rt)$}
	\Let{$\nu$}{$\nu+1$} 
   \EndWhile
   \Ensure{$\balpha^\nu$, $\bB^\nu$.}
  \end{algorithmic}
\end{algorithm}

\begin{figure}[h]
\centering
\includegraphics[width=0.49\textwidth]{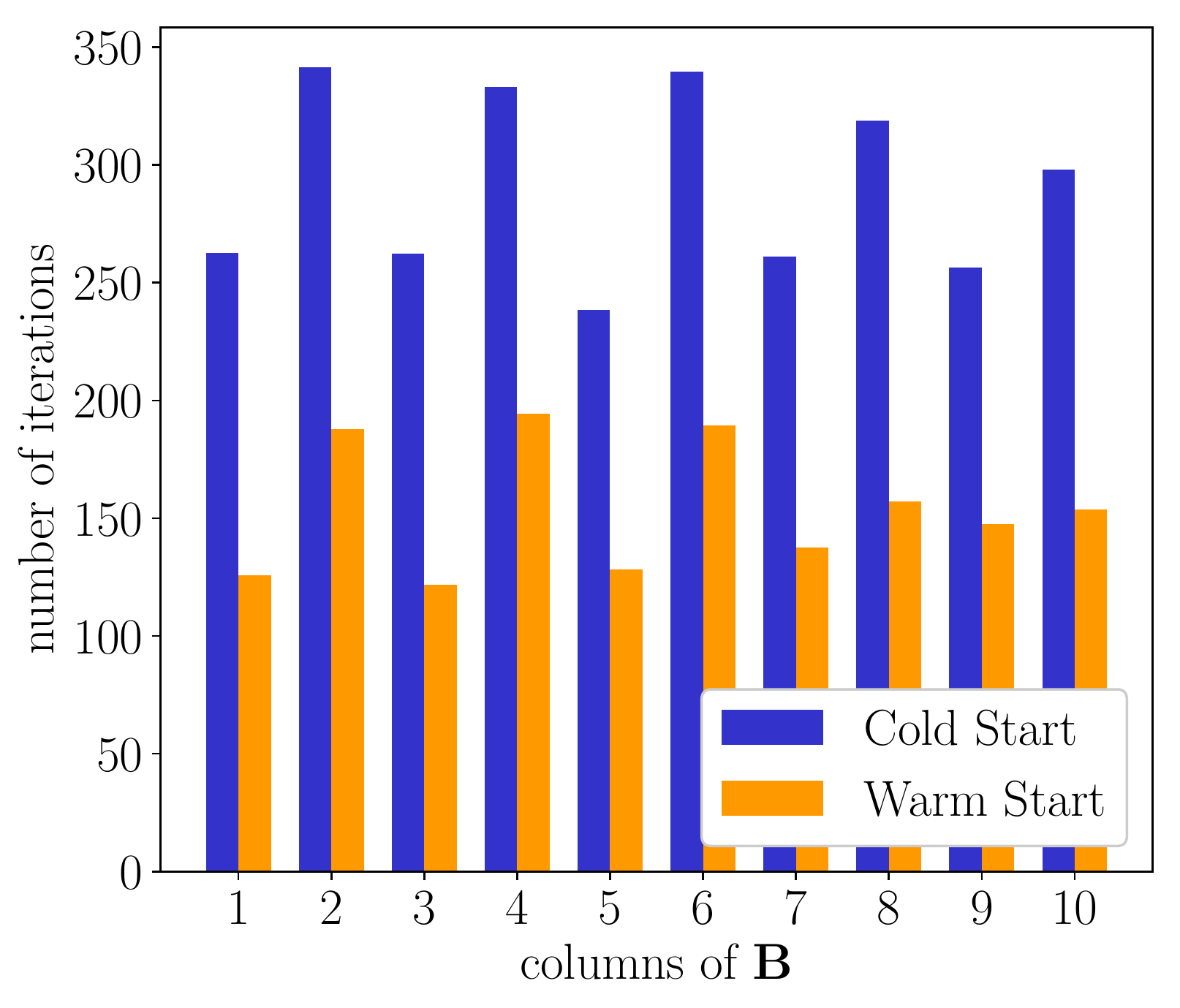}
\caption{\label{fig:bfgs}Average BFGS iterations for each subproblem across the columns.}
\end{figure}

Updating $\bb_j$ can be done efficiently by exploiting the optimized DMD problem structure.  
In particular, BFGS builds a Hessian approximation as it proceeds. 
All of the subproblems for $\bb_j$ share the same $\bPhi(\balpha)$ and, since the columns of $\bB$ contain spatial information, 
neighboring columns are likely to be similar to each other.
After solving subproblem $j$, we use the resulting Hessian approximation to initialize the next subproblem $j+1$. 
Warm starts cut total BFGS iterations in half, see Figure~\ref{fig:bfgs}.

There are different ways to update the weights $\bw$, see line 4 in Algorithms~\ref{alg:BFGS-BFGS} and~\ref{alg:FISTA-BFGS}. 
Define 
\[
\rho_j^\nu = \rho\lt(X_{\cdot j} - \bPhi(\balpha^\nu) \bb_j^{\nu+1}\rt).
\]
The objective with respect to $\bw$ is given by 
\[
\min_{\bw} \sum_{j=1}^n w_j \rho_j^\nu + s(\bw),
\]
where $s$ encodes the weight constraints~\eqref{eq:simplex}. 
The simplest update rule is to set $w_j = 1$ if $\rho_j^\nu$ is one of the $h$ smallest, and $0$ otherwise~\citep{yang2016high}; 
this corresponds to partial minimization in $\bw$ at every step.  
A less aggressive strategy  is to use proximal updates on $\bw$, 
\[
\bw^{\nu+1} = \mbox{prox}_{\eta_\bw s}\lt(\bw^\nu - \eta_\bw \nabla_\bw f(\balpha^\nu, \bB^{\nu+1}, \bw^\nu)\rt)
\] 
where any step size $\eta_\bw > 0$ can be used~\citep{aravkin2016smart}. We use the former simple rule as the default 
in the algorithm. When $h=n$, trimming is turned off, and all weights are identically equal to $1$. 

\begin{algorithm}[h]
\caption{SVRG for DMD}\label{alg:svrg_dmd}
\begin{algorithmic}[1]
\Require{$\balpha^0$, $\bB^0$, $\bw^0$}
\State{Initialize $\nu = 0$, $\zeta_j = \nabla f_j\left(\balpha^0, \bw^0\right)$ for $j = 1,2, \ldots, n$, and $\zeta = \frac{1}{n}\sum_{j=1}^n \zeta_j$}
\While{not converged}
\State{Uniformly sample $I^\nu\subset\{1,2,\ldots,n\}$, such that $\lt|I^\nu\rt| = \tau$}
\State{Sample $J^\nu\in\{0,1\}$, such that $P(J=1) \ll P(J=0)$.}
\For{$j \in I_\nu$}
   		\Let{$\bb^{\nu+1}_j$}{$\displaystyle\arg\min_\bb w_j^\nu \rho(X_{\cdot, j} - \bPhi(\balpha^\nu) \bb) + q(\bb)$}
		\Let{$\zeta_j^+$}{$\nabla \tilde g_j\lt(\balpha^\nu, \bw\rt)$}
	\EndFor
\If{$J = 1$}
\Let{$\bw^{\nu+1}$}{weights update}
\Else
\Let{$\bw^{\nu+1}$}{$\bw^\nu$}
\EndIf
\Let{$\balpha^{\nu+1}$}{$\text{prox}_{\eta_\balpha r}\lt(\balpha^\nu - \eta_\balpha\lt[\frac{1}{\tau}\sum_{j\in I^\nu}\lt(\zeta_j^+ - \zeta_j\rt) + \zeta\rt]\rt)$}
\Let{$\eta_\balpha$}{step size update}
\Let{$\zeta_j$}{$\zeta_j^+$  for $j\in I^\nu$}
\Let{$\zeta$}{$\frac{1}{n}\sum_{j=1}^n \zeta_j$}
\Let{$\nu$}{$\nu+1$}
\EndWhile
\Ensure{$\balpha^\nu$, $\bB^\nu$}
\end{algorithmic}
\end{algorithm}
\subsection{A scalable stochastic method}
In DMD applications, $n$ represents the number of spatial variables, and is often much larger than either 
$m$ or $k$. Therefore, step 2 of Algorithms~\ref{alg:BFGS-BFGS} and~\ref{alg:FISTA-BFGS} is a computational bottleneck. 
We use stochastic methods to scale the approach. The basic idea is to partially minimize  
over a random sample of the columns of $\bB$; the resulting (scaled) gradient is an unbiased estimate of $\nabla_{\balpha}\tilde f$. 
More precisely, define 
\begin{align*}
\bb_j(\balpha, \bw) &= \argmin_\bb w_j\rho\lt(X_{\cdot j} - \bPhi(\balpha) \bb\rt) + q(\bb),\\
\tilde g_j(\balpha, \bw) &= w_j\rho(X_{\cdot j} - \bPhi(\balpha) \bb_j(\balpha, \bw)) + q(\bb_j(\balpha, \bw)).
\end{align*}
Then we have
\[
\tilde f(\balpha, \bw) = \sum_{j=1}^n \tilde g_j(\balpha, \bw) + r(\balpha) + s(\bw).
\]
This is a classical setting for stochastic methods. In each iteration, we can use a subset of $\tilde g_j$ to calculate the approximate gradient for the smooth part of $\tilde f$ in order to reduce the computational burden. Here we use SVRG~\citep{johnson2013accelerating} as our stochastic solver for the outer problem; the full details are given in Algorithm~\ref{alg:svrg_dmd}. Note that this
stochastic approach is an alternative to using a cost reduction
based on projecting onto SVD modes \citep{askham2017variable} or using an optimized 
but fixed subsampling of the columns \citep{gueniat2015a}. With the method
of Algorithm~\ref{alg:svrg_dmd}, 
none of the data is discarded or filtered by the cost reduction procedure.

\begin{figure}[h]
\centering
\includegraphics[width=0.6\textwidth]{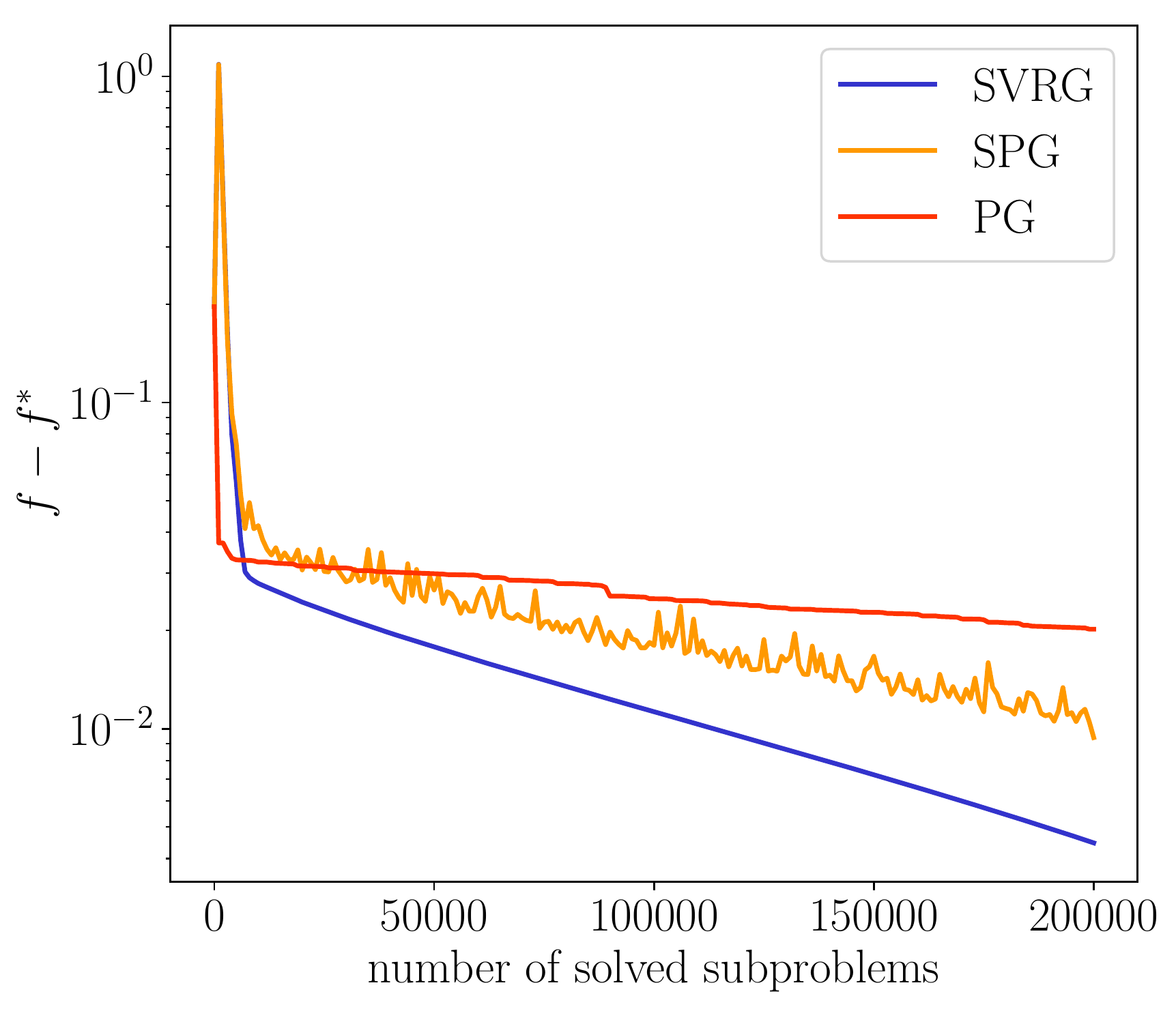}
\caption{\label{fig:compare}Compare performance of SVRG, Stochastic Proximal Gradient (SPG) method and Proximal Gradient (PG) method over the same data set.}
\end{figure}

In Figure~\ref{fig:compare}, we solve a problem with dimension $m=512$ and $n=1000$.
A diminishing step size scheme is used, taking 
\[
\eta_\balpha^\nu = \frac{\eta_\balpha^0}{\mbox{floor}(\nu/K) + 1}, 
\]
with $\eta^0_\balpha = 10^{-7}$ and $ K = 500$ for the result in Figure~\ref{fig:compare}.
Comparing the algorithms according to total $\bb_j$ subproblems, we see that SVRG decreases faster and is less noisy than
the Stochastic Proximal Gradient (SPG) method\footnote{It is important to note that SPG has no convergence theory, 
while SVRG is guaranteed to converge.   In practice SPG works well so we include it in the comparison.}. Proximal Gradient (PG) decreases quickly 
in the beginning, but is soon overtaken by stochastic methods. SVRG gives a significant improvement over SPG. 

The trimming weights $\bw$ rely on global information; that is, the best $h$ residuals are easily selected after all of the residuals 
have been calculated. This is why the weights update (lines 8-11 of Algorithm~\ref{alg:svrg_dmd}) is done rarely.  
For detailed analysis of stochastic algorithms with trimming, see~\cite{aravkin2016smart}.




%

\section{Synthetic examples}
\label{sec:synthetic}
In this section, we examine the performance of the
robust DMD on a pair of synthetic test cases
with known solution. These examples are drawn from the
additive noise study of \citet{dawson2016}.

\subsection{A simple periodic example}

\begin{figure}
  \centering
  \includegraphics[width=0.49\textwidth]{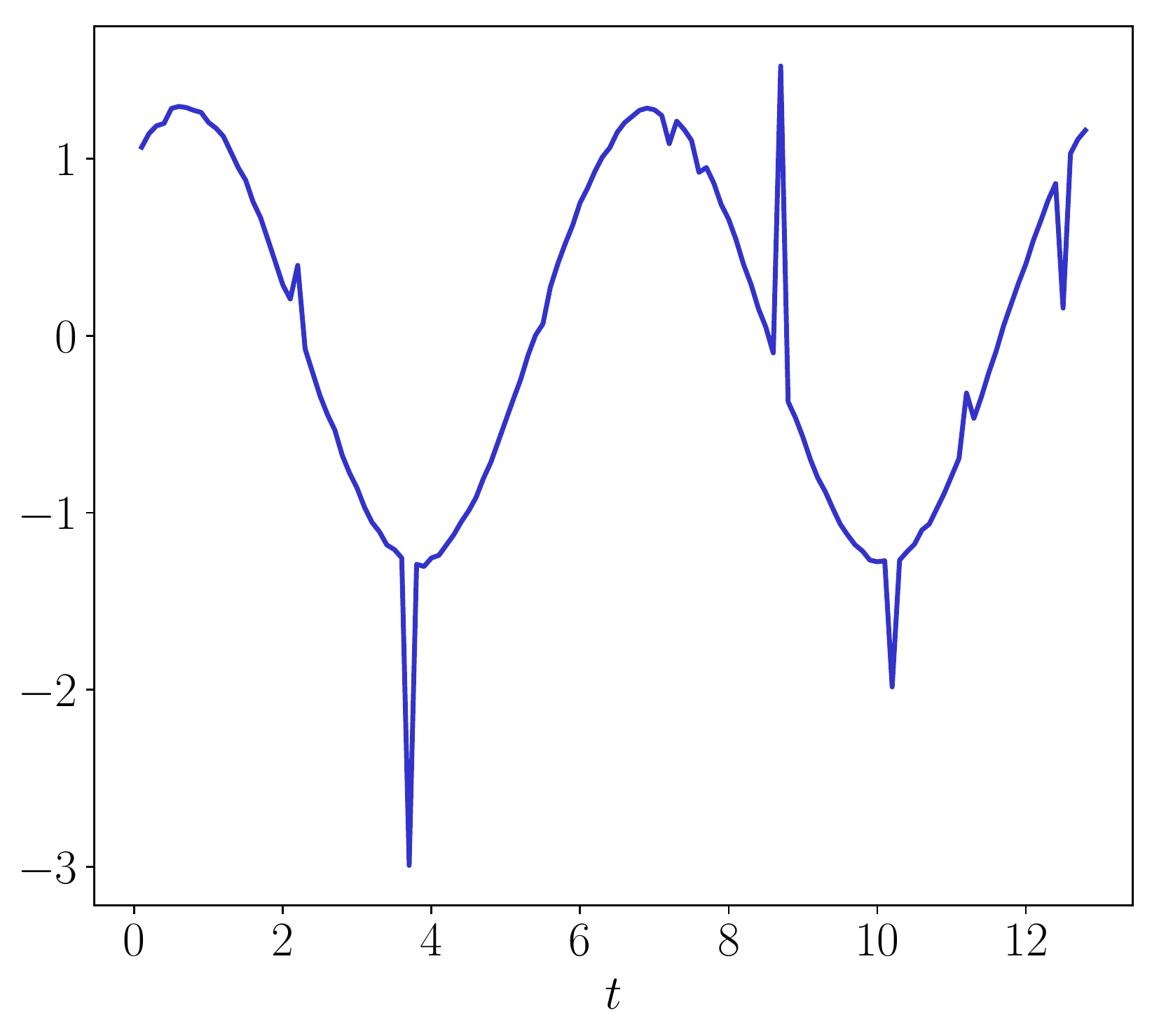}
  \includegraphics[width=0.49\textwidth]{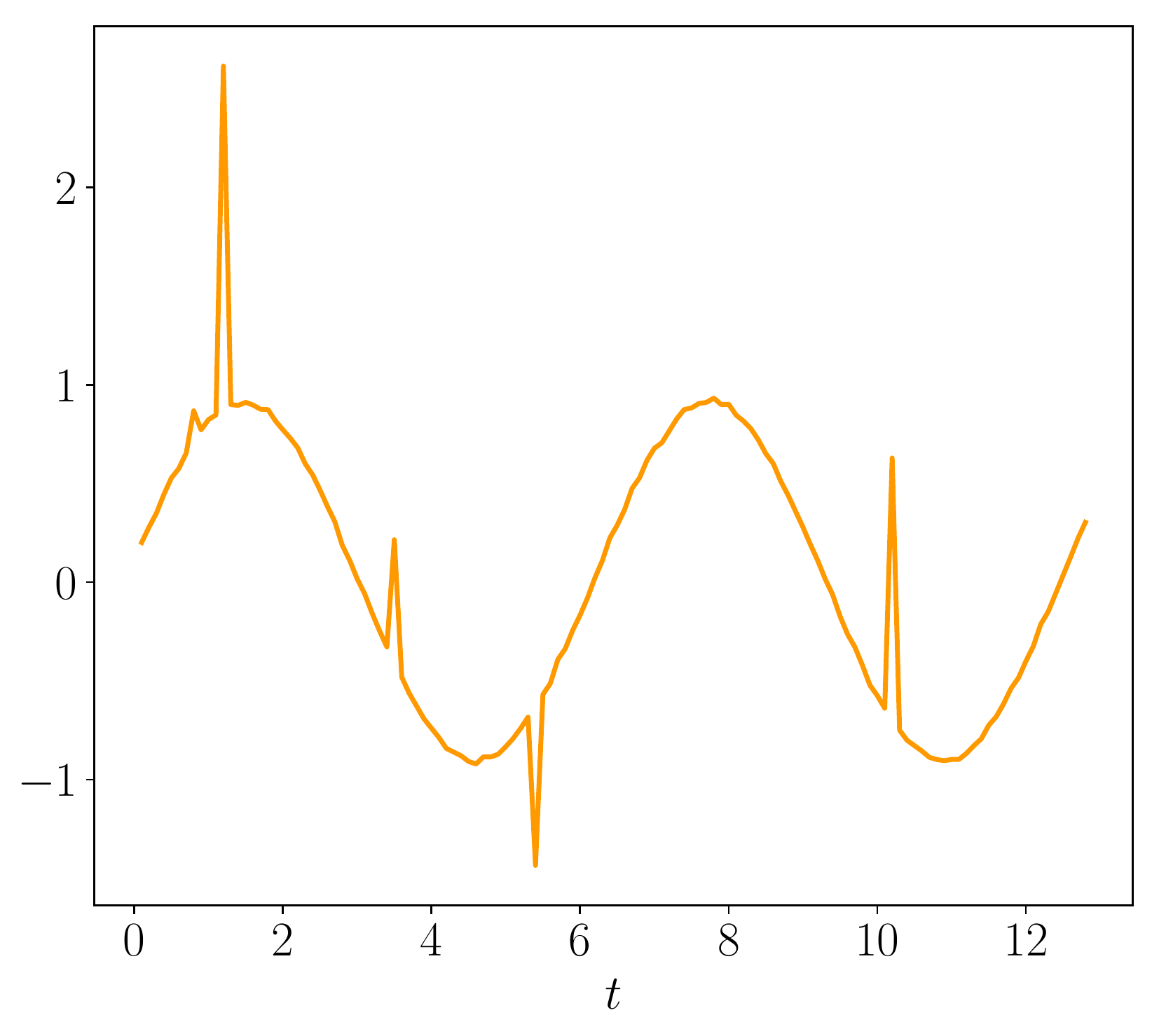}
  \caption{Sample time series of $x_1(t)$ and $x_2(t)$ for the
    simple periodic example, with background noise of size
    $\sigma = 10^{-2}$ and spikes of size $\mu = 1$ added at
    $p = 5\%$ of the snapshots for each channel.}
  \label{fig:test1_sample}
\end{figure}

Let $\bx(t)$ be the solution of a two dimensional
linear system with the following dynamics

\begin{equation} \label{eq:test1}
  \dot{\bx} = \begin{pmatrix} 1 & -2 \\ 1 & -1 
  \end{pmatrix} \bx \; .
\end{equation}
We use the initial condition $\bx(0) = (1,0.1)^\intercal$
and take snapshots

\begin{equation}
  \bx_j = \bx(j \Delta t)
  + \sigma {\mathbf g}_j + \mu {\mathbf s}_j \, , \nonumber
\end{equation}
where $\Delta t = 0.1$, $\sigma$ and $\mu$ are prescribed
noise levels, ${\mathbf g}_j$ is a vector whose entries are
drawn from a standard normal distribution, and ${\mathbf s}_j$
is a vector whose entries are the product of a
Bernoulli trial with small expectation
$p$ and a standard normal (corresponding to sparse
noise). The snapshots are therefore corrupted with a
base level of noise $\sigma$ and sparse ``spikes'' of
size $\mu$ with firing rate $p$. A sample time series
for this example can be found in Figure~\ref{fig:test1_sample}.

\begin{figure}
  \centering
  \includegraphics[width=0.5\textwidth]{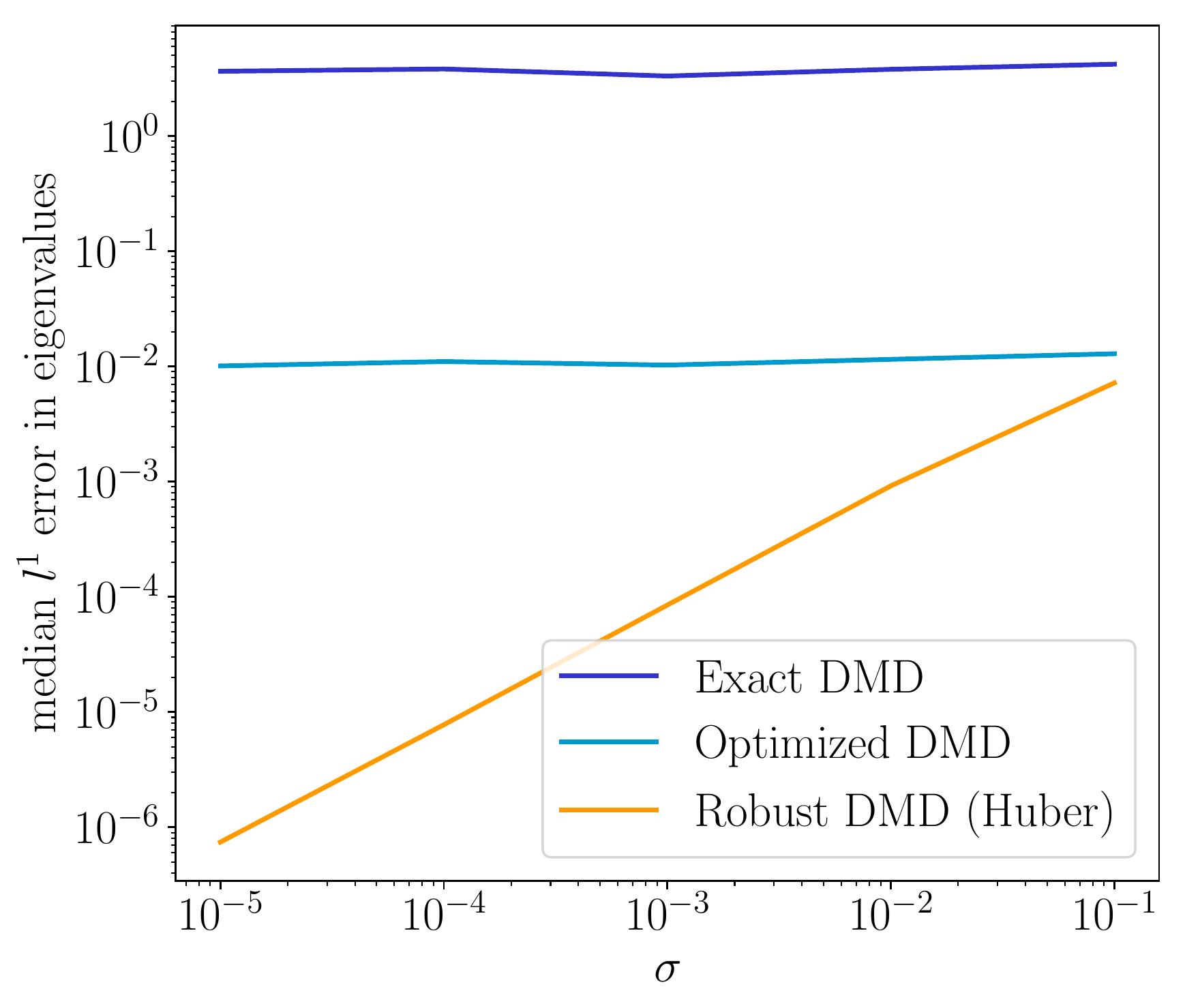}
  \caption{Median error in the computed eigenvalues
    over 200 runs. The background noise $\sigma$ varies
    while the size of the spikes is fixed at
    $\mu = 1$ and the firing rate is fixed at $p=5\%$.}
  \label{fig:test1_robdmdsparsenoise}
\end{figure}

The $k=2$ eigenvalues of the system matrix in~\eqref{eq:test1}
are $\pm \bi$, corresponding to sinusoidal dynamics in
time. In Figure~\ref{fig:test1_robdmdsparsenoise}, we plot
the median (over 200 random trials) of the $l^1$-norm error in the
approximations of these eigenvalues using three different
methods: the exact DMD of \citet{tu2014dynamic}; the optimized DMD
as defined in \eqref{eq:optdmd};
and the robust DMD as defined in \eqref{eq:problemform},
with $\rho$ the Huber norm and $h = n = 2$ (no trimming).
Each trial consists of the first 128 snapshots with additive
noise.
The level of the background noise, $\sigma$, varies
over the experiments and the size and firing rate of
the spikes are fixed at $\mu=1$ and $p=5\%$, respectively.
We set the Huber parameter using knowledge of the problem
set-up, i.e. $\kappa = 5 \sigma$, but in a real-data
setting this parameter would have to be estimated or
chosen adaptively. While the optimized DMD improves over
the exact DMD, the error does not decrease as the level
of the background noise decreases. We therefore see
the effect of the sparse outliers using the optimized DMD.
For the robust
formulation, on the other hand, the accuracy of the eigenvalues
is determined by the level of the background noise,
so that the outliers are not biasing the computed
eigenvalues.

\subsection{An example with hidden dynamics}

\begin{figure}
  \centering
  \begin{subfigure}[b]{0.49\textwidth}
    \includegraphics[width=\textwidth]{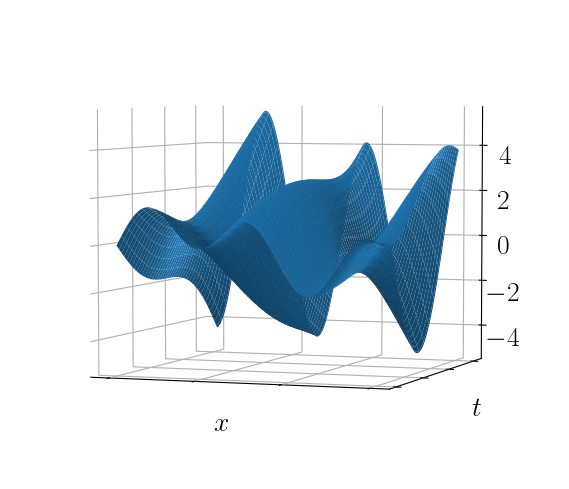}
    \caption{clean data}
    \label{fig:test2samp_xclean}
  \end{subfigure}
  \begin{subfigure}[b]{0.49\textwidth}
    \includegraphics[width=\textwidth]{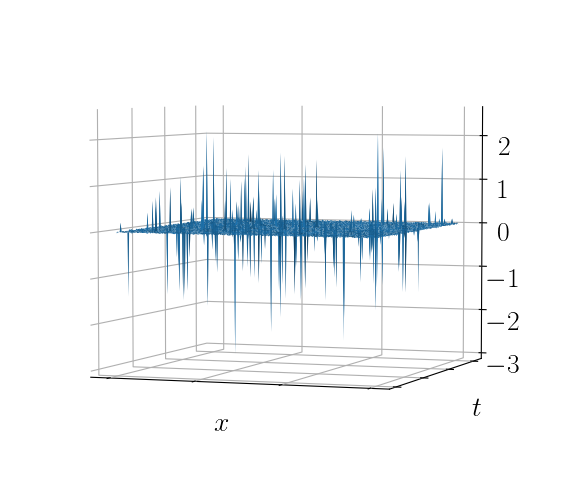}
    \caption{``sparse noise''}
    \label{fig:test2samp_spikes}
  \end{subfigure}
  \begin{subfigure}[b]{0.49\textwidth}
    \includegraphics[width=\textwidth]{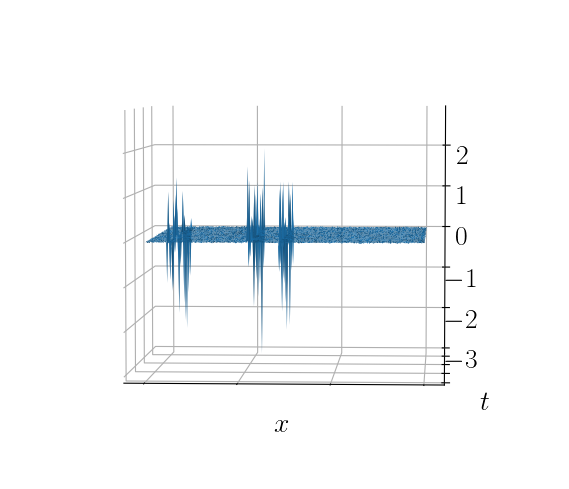}
    \caption{``broken sensor''}
    \label{fig:test2samp_col}
  \end{subfigure}
  \begin{subfigure}[b]{0.49\textwidth}
    \includegraphics[width=\textwidth]{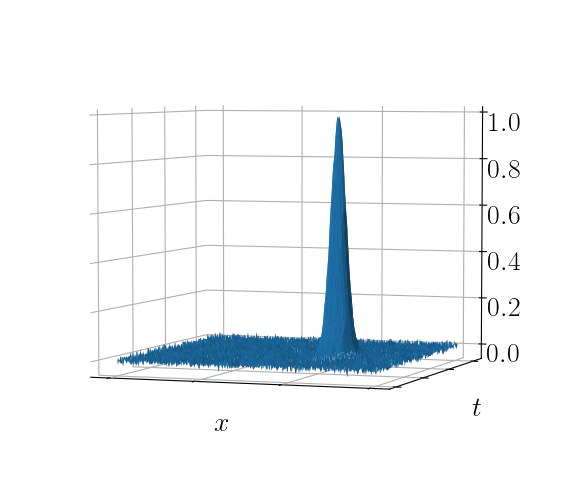}
    \caption{``bump''}
    \label{fig:test2samp_bump}
  \end{subfigure}
  \caption{A surface plot of the data for the hidden dynamics
    example and surface plots of a sample of each type of 
    noise we consider.}
  \label{fig:test2_samp}
\end{figure}

In the case that a signal contains some rapidly decaying
components it can be more difficult to identify  the
dynamics, particularly in the presence of sensor noise
\citep[see][]{dawson2016}. We consider a signal composed of two
sinusoidal forms which are translating, with one growing
and one decaying, i.e.

\begin{equation}
  x(y,t) = \sin( k_1y - \omega_1 t) e^{\gamma_1 t} +
  \sin( k_2y - \omega_2 t) e^{\gamma_2 t} \; ,
\end{equation}
where $k_1 = 1$, $\omega_1 = 1$, $\gamma_1 = 1$, $k_2 = 0.4$,
$\omega_2 = 3.7$, and $\gamma_2 = -0.2$ (following
settings used by \cite{dawson2016}). This signal
has $k=4$ continuous time eigenvalues given
by $\gamma_1 \pm i \omega_1$ and $\gamma_2 \pm i \omega_2$.
We set the domain of $y$ to be $[0,15]$ and use
300 equispaced points, $y_j$, to discretize. For the
time domain, we set $\Delta t = \pi/(2^8-2)$
so that the number of snapshots we use,
$m = 2^7$, covers $[0,\pi/2]$. We denote the 
vector of discrete values $x(y_j,t)$ by $\bx(t)$.
See Figure~\ref{fig:test2samp_xclean} for a surface plot
of this data.

We consider three different types of perturbations of 
the data. The first perturbation adds background noise
and spikes, as in the previous example, i.e. the snapshots
are given by

\begin{equation}
  \bx^{(1)}_j = \bx(j \Delta t)
  + \sigma {\mathbf g}_j + \mu {\mathbf s}_j \, , \nonumber
\end{equation}
where $\sigma$ and $\mu$ are prescribed
noise levels, ${\mathbf g}_j$ is a vector whose entries are
drawn from a standard normal distribution, and ${\mathbf s}_j$
is a vector whose entries are the product of a
Bernoulli trial with small expectation
$p$ and a standard normal. See Figure~\ref{fig:test2samp_spikes}
for a sample plot of this ``sparse noise'' pattern. The second 
perturbation we consider adds background noise and spikes which
are confined to specific entries of $\bx_j$, i.e. 

\begin{equation}
  \bx^{(2)}_j = \bx(j \Delta t)
  + \sigma {\mathbf g}_j + \mu \tilde{{\mathbf s}}_j \, , \nonumber
\end{equation}
where ${\mathbf g}_j$, $\sigma$, and $\mu$ are as above 
and the $\tilde{\mathbf s}_j$ are sparse vectors
which have the same sparsity pattern for all $j$ and
nonzero entries drawn from a standard normal distribution 
(this corresponds to having a few broken sensors recording 
the data). We plot a sample of this ``broken sensor'' noise
pattern in Figure~\ref{fig:test2samp_col}. The third perturbation 
we consider adds background noise and a localized bump to the data, 
i.e. 

\begin{equation}
  \left [ \bx^{(3)}_j \right]_i = x(y_i, j \Delta t)
  + \sigma \mathcal{N}(0,1) 
  + A \exp \left ( -\left ( \frac{y_b-y_i}{w\Delta y} \right)^2 
  - \left ( \frac{t_b-j\Delta t}{w\Delta t} \right ) ^2 \right ) \, 
  , \nonumber
\end{equation}
where $\sigma$ is as above, $\mathcal{N}(0,1)$ denotes a number
drawn from the standard normal distribution, $A$ determines 
the maximum height of the bump, $w$ determines the ``width'' 
of the bump, and $y_b$ and $t_b$ determine the center of the bump 
in space and time (this corresponds to having some non-exponential 
dynamics in the data). In Figure~\ref{fig:test2samp_bump}, we 
plot a sample of this ``bump'' noise pattern.

\begin{figure}
  \centering
  \begin{subfigure}[b]{0.49\textwidth}
    \includegraphics[width=\textwidth]{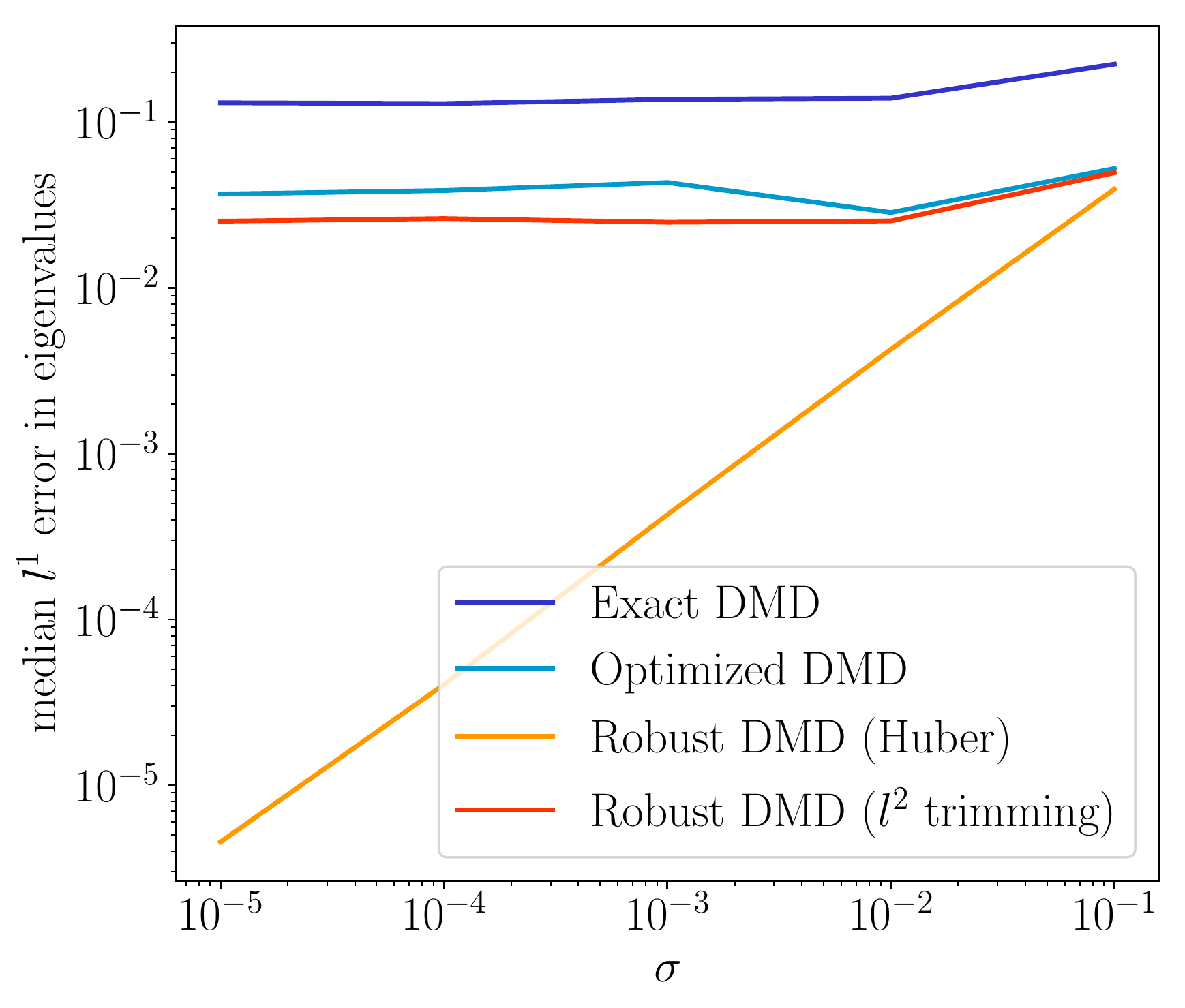}
    \caption{``sparse noise''}
    \label{fig:test2_errs_spikes}
  \end{subfigure}
  \begin{subfigure}[b]{0.49\textwidth}
    \includegraphics[width=\textwidth]{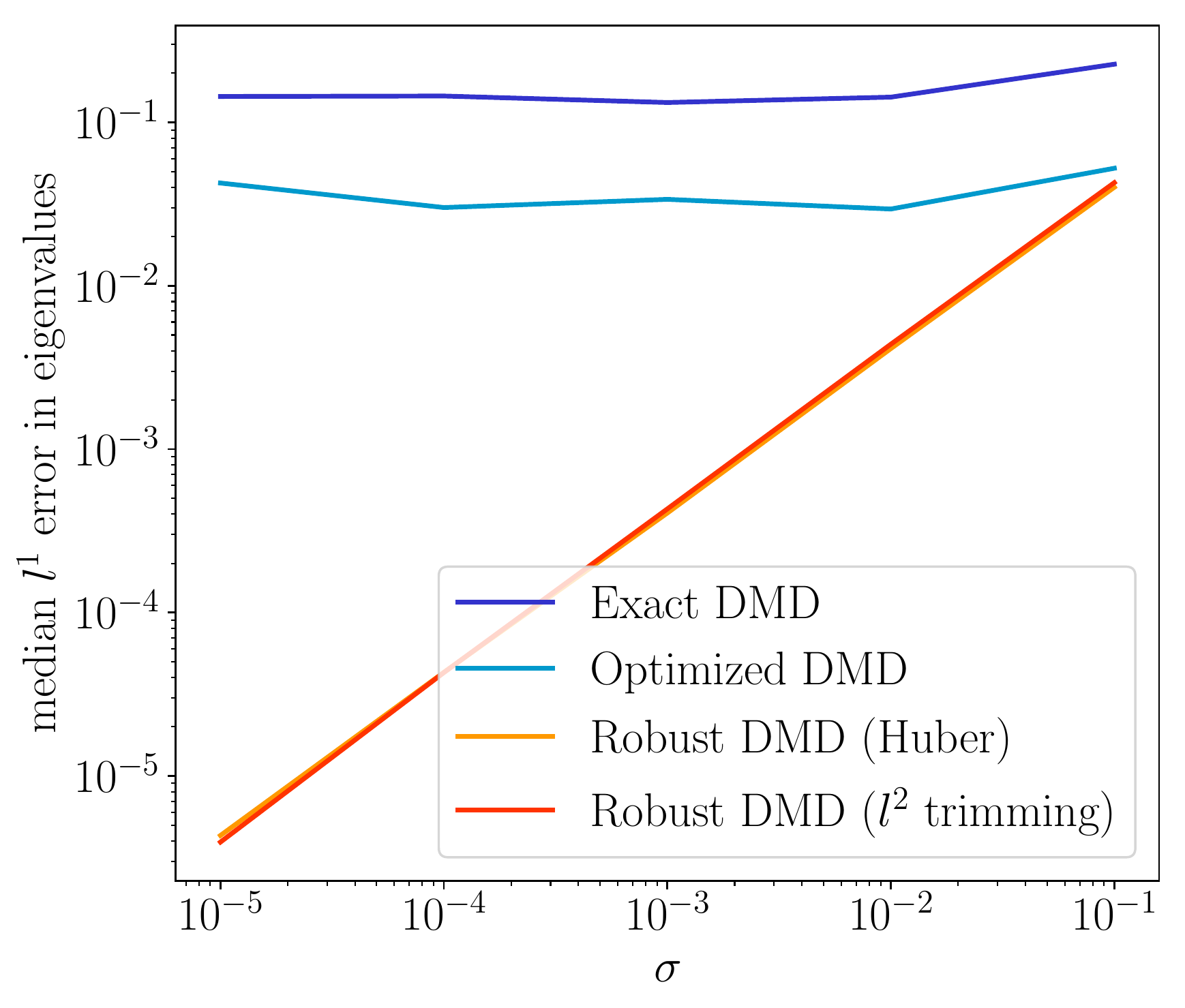}
    \caption{``broken sensor''}
    \label{fig:test2_errs_col}
  \end{subfigure}
  \begin{subfigure}[b]{0.49\textwidth}
    \includegraphics[width=\textwidth]{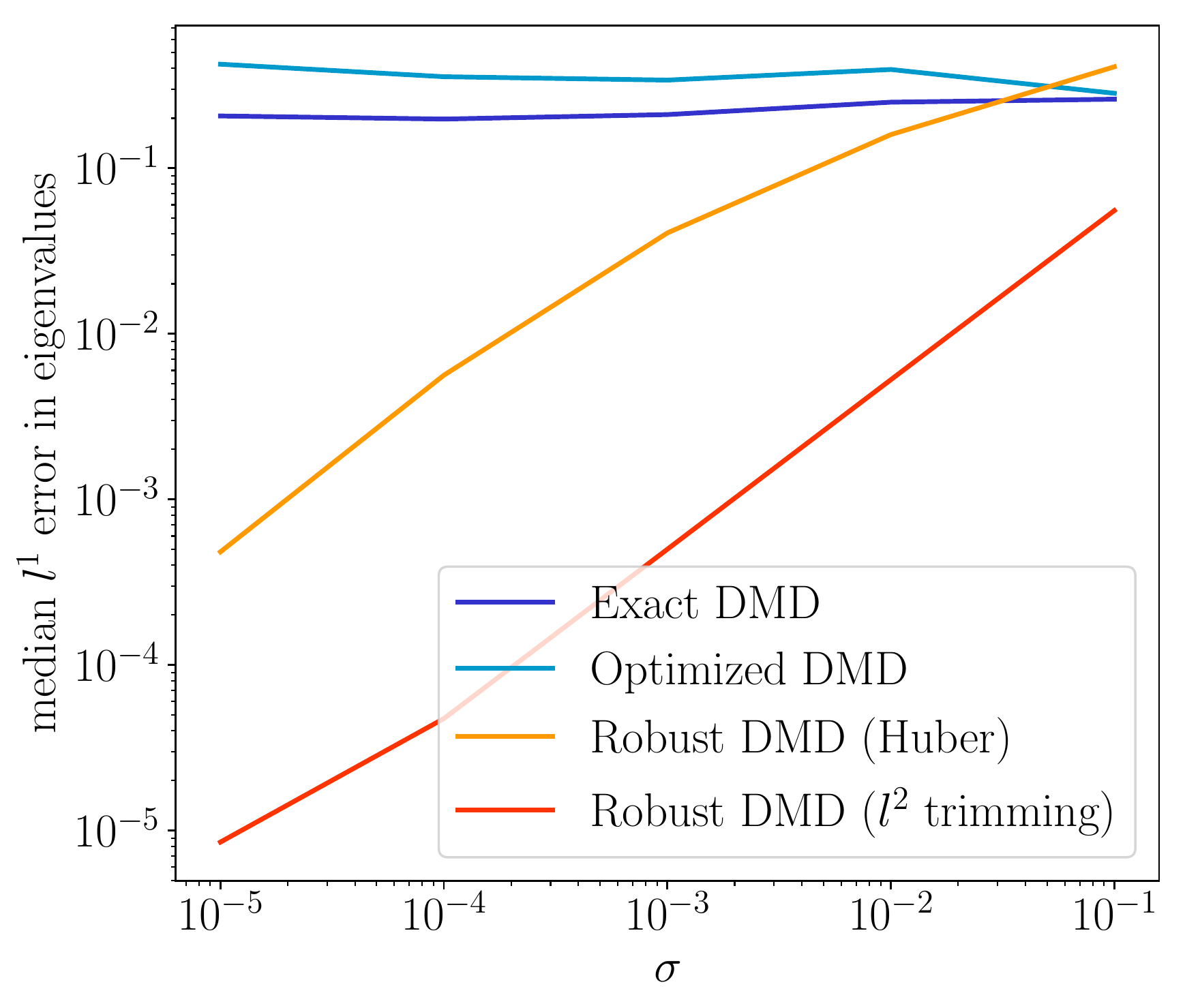}
    \caption{``bump''}
    \label{fig:test2_errs_bump}
  \end{subfigure}
  \caption{Median error in the computed eigenvalues
    over 200 runs. The background noise $\sigma$ varies
    while the size of the spikes is fixed at
    $\mu = 1$ and the firing rate is fixed at $p=5\%$
    for the ``sparse noise'' and ``broken sensor'' 
    examples and the height is fixed at $A=1$ and 
    the width at $w=10$ for the ``bump'' example.}
  \label{fig:test2_errs}
\end{figure}

In Figure~\ref{fig:test2_errs}, we plot
the median (over 200 random trials) of the $l^1$-norm error in the
approximations of the eigenvalues using four different
methods: the exact DMD of \citet{tu2014dynamic}; the optimized DMD
as defined in \eqref{eq:optdmd};
the robust DMD as defined in \eqref{eq:problemform},
with $\rho$ the Huber norm and $h = n = 300$ (no trimming);
and the robust DMD with $\rho$ the standard Frobenius norm and 
$h = 0.8n = 240$ (trimming). Each trial consists of the first 128 
snapshots with additive noise.
The level of the background noise, $\sigma$, varies
over the experiments. For the ``sparse noise'' and ``broken sensor''
snapshots, the size of the spikes is fixed at
$\mu = 1$ and the density is fixed at $p=5\%$, i.e. $5\%$ 
of the entries are corrupted for the ``sparse noise'' example 
and $5\%$ of the sensors are corrupted for the ``broken sensor''
example.
For the ``bump'' snapshots, the height of the bump is fixed at $A=1$ 
and the width  at $w=10$.
We set the Huber parameter using knowledge of the problem
set-up, i.e. $\kappa = 5 \sigma$, but in a real-data
setting this parameter would have to be estimated or
chosen adaptively. 

With sparse noise, as in 
Figure~\ref{fig:test2_errs_spikes}, the results for the exact DMD, 
optimized DMD, and Huber norm-based robust DMD are consistent 
with the simple periodic example. The Huber norm formulation is
the only one which is able to take advantage of the lower levels
of background noise. The trimming formulation provides very 
little advantage for this example, as any sensor can be affected
by the outliers. In contrast, we see that the trimming formulation
is able to perform as well as the Huber formulation for the 
broken sensor example  (see Figure~\ref{fig:test2_errs_col}), 
as the algorithm is able to adaptively remove the broken sensors 
from the data. In Figure~\ref{fig:test2_errs_bump}, we plot the 
results for the bump data, which display some interesting behavior.
The optimized DMD actually performs worse than the exact DMD, which
is attributable to over-fitting. For all but the highest background
noise level, the Huber and trimming formulations show a significant 
advantage over the optimized DMD and exact DMD, with the trimming 
formulation performing the best. The trimming formulation therefore
presents an attractive solution for data with unknown, localized 
deviations from the exponential basis of the DMD, especially given
that the inner problem for trimming with the Frobenius penalty 
can be solved rapidly. Of course, trimming can be combined with
a Huber (or other robust) penalty for increased robustness to 
outliers.

\section{Real data examples}
\label{sec:real}
These examples will be added to a future draft.
\subsection{Atmospheric chemistry data}

\subsection{Fluid simulation data}

\section{Conclusion and future directions}
\label{sec:conclusion}

We have presented an optimization framework
and a suite of numerical algorithms for computing the
dynamic mode decomposition with robust penalties and
parameter constraints. This framework allows for improved
performance of the DMD in a number of settings, as
borne out by synthetic and real data experiments.
In the presence of sparse noise or non-exponential
structure, the use of robust penalties significantly
decreases the bias in the computed eigenvalues. When
using the DMD to perform future state prediction, adding
the constraint that the eigenvalues lie in the left
half-plane increases the stability of the extrapolation.
The algorithms presented are capable of solving small
to medium-sized problems in seconds
on a laptop and scale well
to higher-dimensional problems due to their intrinsic
parallelism and the efficiency of the SVRG approach.
In contrast with previous approaches, the SVRG increases
efficiency without throwing out data or incidentally
filtering it. We believe that the framework and
algorithms presented here will enable practitioners of
the DMD to tackle larger, noisier, and more complex
data sets than previously possible. The authors commit
to releasing the software used for these calculations
as an open-source package in the Julia language
\citep{bezanson2012julia}.

The present work can be extended in a number of
ways. Because the inner solve completely decouples
over the columns of ${\bf X}$ and ${\bf B}$, the
algorithms presented above immediately generalize to
data-sets with missing entries and even data which
are collected asynchronously across sensors.
While the global nature of an optimized DMD
fit has advantages in terms of the quality of the
recovered eigenvalues, it implicitly rules out
process noise. However, including process noise or
a known forcing term would be useful in many
applications. Incorporating such terms into
this optimization framework is ongoing work and results
will be reported at a later date. We also note
that much of the above applies to dimensionality
reduction using any parameterized family of time
dynamics, not just exponentials. For such an
application, many of the algorithms above could be
easily adapted, so long as gradient formulas are
available.


\acks{T. Askham and J. N. Kutz acknowledge support from 
the Air Force Office of Scientific Research (FA9550-15-1-0385).  
J. N. Kutz also acknowledges support from the Defense Advanced 
Research Projects Agency (DARPA contract HR0011-16-C-0016).
The work of A. Aravkin and P. Zheng was supported by the Washington 
Research Foundation Data Science Professorship. 
}

\bibliography{refs}

\end{document}